\documentclass[a4paper,12pt,twoside]{article}
\usepackage{graphicx}\usepackage{a4wide}
\usepackage{amssymb,amsmath,amsthm}
\usepackage[only,llbracket,rrbracket,interleave]{stmaryrd}
\usepackage[pdftex,colorlinks,bookmarksnumbered,bookmarksopen,citecolor=red,urlcolor=red]{hyperref}
\usepackage{cite}

\usepackage[overload]{empheq}
\newcommand{\bm}[1]{\text{\boldmath $#1$\unboldmath}}

\newcommand{\vect}[1]{\mathbf{#1}}
\newcommand{\mat}[1]{\mathbf{#1}}

\newcommand{\RR}{\mathbb{R}}

\newcommand{\sobo}[1][1]{\ensuremath{\mathcal{H}^{#1}}}
\newcommand{\eltwo}{\ensuremath{\mathcal{L}_2}}
\newcommand{\trial} {\ensuremath{\mathcal{V}}}
\newcommand{\test} {\sobo_{\Gamma_D}}
\newcommand{\nsd}  {\ensuremath{\texttt{n}_{\texttt{sd}}}}

\newcommand{\npar}    {\ensuremath{\texttt{n}_{\texttt{pa}}}}

\newcommand{\bI}{\bm{\mathcal{I}}}
\newcommand{\I}{\mathcal{I}}

\newcommand{\bn}{\bm{n}}
\newcommand{\bx}{\bm{x}}
\newcommand{\grad}{\bm{\nabla}}
\newcommand{\bu}{\bm{u}}
\newcommand{\bv}{\bm{v}}
\newcommand{\bw}{\bm{w}}
\newcommand{\bmu}{\bm{\mu}}

\newcommand{\upgd}{\bu_{_{\texttt{PGD}}}}
\newcommand{\ppgd}{p_{_{\texttt{PGD}}}}
\newcommand{\uupgd}{\tilde{\bu}_{_{\texttt{PGD}}}}
\newcommand{\pppgd}{\tilde{p}_{_{\texttt{PGD}}}}

\newcommand{\Fu}{\bm{F}_{\!\!\bu}}
\newcommand{\Lu}{L_{\bu}}
\newcommand{\Fp}{F_{\! p}}
\newcommand{\Lp}{L_{p}}
\newcommand{\de}{\delta\!}

\newcommand{\SZ}[1]{{\color{black} #1\color{black}}}
\newtheorem{remark}{Remark}

\begin{document}

\title{Generalized parametric solutions in Stokes flow}

\author{Pedro D\'{\i}ez  \and Sergio Zlotnik \and Antonio Huerta\footnote{Corresponding author: A. Huerta, Laboratori de C\`alcul Num\`eric (LaC\`aN). ETS de Ingenieros de Caminos, Universitat Polit\`ecnica de Catalunya, Jordi Girona 1, 08034 Barcelona, Spain.}}
 \date{Laboratori de C\`alcul Num\`eric (LaC\`aN).\\ 
 ETS de Ingenieros de Caminos, Canales y Puertos, \\
 Universitat Polit\`ecnica de Catalunya\textperiodcentered BarcelonaTech, Barcelona, Spain.\\}


\maketitle
\begin{abstract}
Design optimization and uncertainty quantification, among other applications of industrial interest, require fast or multiple queries of some parametric model. The Proper Generalized Decomposition (PGD) provides a separable solution, \SZ{a \emph{computational vademecum}} explicitly dependent on the parameters, efficiently computed with a greedy algorithm combined with an alternated directions scheme and compactly stored. This strategy has been successfully employed in many problems in computational mechanics. The application to problems with saddle point structure raises some difficulties requiring further attention. This article proposes a PGD formulation of the Stokes problem. Various possibilities of the separated forms of the PGD solutions are discussed and analyzed, selecting the more viable option. The efficacy of the proposed methodology is demonstrated in numerical examples for both Stokes and Brinkman models.

\end{abstract}


\section{Introduction}

Standard discretization techniques in computational mechanics have reached an amazing level of maturity and efficiency. Nonetheless, the systematic exploration of parametric solutions arising from optimization (where the best choice for the parameters is unknown) or uncertainty quantification (where the parameters have stochastic features) is often computationally unaffordable. The Proper Generalized Decomposition (PGD), among other Reduced Order Models, provides a generalized solution with an explicit parametric dependence. This compact expression containing the analytical dependence on the free parameters is also known as \emph{computational vademecum} and allows an expedited exploration of the parametric space, with the computational cost of a simple interpolation, i.e. post-processing. 

The PGD has been successfully employed in different problems in the broad field of mathematical and computational modelling. Essentially, PGD consists in finding a separable approximation, that is a sum of terms, each of them being a product of modal functions depending on one of the parameters. This approximation is usually computed with a greedy algorithm (obtaining the terms sequentially) and, for each term, an alternated directions iterative scheme is, in general, employed to find the different parametric modes.

PGD was already used in the framework of Stokes and Navier--Stokes governed problems to obtain separated solutions in terms of the different spatial dimensions for Cartesian domains  \cite{Dumon:11,Dumon:13,Dumon13b} and also for space-time separation \cite{Aghighi:13,Leblond:14}. PGD for space-space separation and space-time separation is significantly increasing the computational efficiency in solving problem with complex flow patterns in simple cartesian domains. 

\SZ{Here, the focus is on solving parametric problems in complex domains with arbitrary geometries  (not assumed to be Cartesian). Therefore the space coordinates are treated together and separated from the different (independent) parametric dimensions.}

This paper aims at analyzing the application of PGD to problems with saddle point structure, taking the Stokes problem as one of the simplest. In particular, special attention is paid to the selection of the form of parametric separation in hybrid formulations. In other words, in a velocity-pressure formulation, the question is: must the parametric modes be independent for velocity and pressure, or just one for both?

Thus, the remainder of the paper is structured as follows. The parametric version of the Stokes problem is stated in Section \ref{sec:prob}. Section \ref{sec:PGD} describes the general formulation of PGD in this case and analyzes the possible alternatives for the parametric separation, concluding that the same parametric mode must be used for both velocity and pressure. Section \ref{sec:exa} presents examples demonstrating the viability of the devised approach.

\section{The parameterized Stokes flow}\label{sec:prob}

The strong form of the Stokes problem can be written as
\begin{equation} \label{eq:Stokes}
 \left\{\begin{aligned}
 -\grad\cdot\nu\grad\bu + \grad p &= \bm{b} &&\text{in $\Omega$}\\
                               \grad\cdot\bu &= 0         &&\text{in $\Omega$}\\
                                               \bu &= \bu_D &&\text{on $\Gamma_D$}\\
        -p\bn+\nu\bn\cdot\grad\bu &= \bm{t}    &&\text{on $\partial\Omega\setminus\Gamma_D$}.
 \end{aligned}\right.
\end{equation}
%

User-prescribed data are the computational domain $\Omega\subset\RR^{\nsd}$ ($\nsd$ being the number of spatial dimensions) whose boundary $\partial\Omega$ is partitioned into Dirichlet, $\Gamma_D$, and Neumann frontiers, the body forces $\bm{s}$, the Dirichlet, $ \bu_D$, and Neumann, $\bm{t}$, boundary conditions, and the kinematic viscosity $\nu$.

Any of these user-prescribed data could be a function of a set of parameters $\bmu\in\bI\subset\mathbb{R}^{\npar}$ (with $\npar$ number of parameters). Those affecting the right-hand-side of the resulting equations (viz. $\bm{s}$, $ \bu_D$ and $\bm{t}$) are easy to handle. On the contrary those affecting the differential operator (viz.\ viscosity or domain) cannot be treated trivially. The set $\bI\subset\mathbb{R}^{\npar}$, which characterizes the admissible range for parameters $\bmu$, can be defined as the cartesian product of the range for each parameter, namely,  $\bI := \I_1\times\I_2\times\dotsb\times\I_{\npar}$ with $\mu_i\in\I_i$ for $i=1,\dotsc , \npar$. 

This can be interpreted as taking $\bmu$ as additional independent variables (or parametric coordinates) instead of problem parameters. Hence, the unknown velocity-pressure pair $(\bu,p)$ can be seen as functions in a larger dimensional space and can be expressed as $\bu(\bx,\bmu)$ and $p(\bx , \bmu )$ with  $(\bx , \bmu ) \in\Omega\times\bI$. 

Consequently, formally $\bu$ and $p$ lie in tensor product spaces, namely
\begin{equation}
 \left\{\begin{aligned}
 \bu&\in\bm{\trial}
 \SZ{        :=[\trial \otimes \eltwo(\I_1) \otimes \eltwo(\I_2) \otimes\dotsb\otimes \eltwo(\I_{\npar})]^{\nsd} } \text{, and}\\
    p&\in\eltwo(\Omega\times\bI) 
         = \eltwo(\Omega)\otimes\eltwo(\I_1)\otimes\eltwo(\I_2)\otimes\dotsb\otimes\eltwo(\I_{\npar}),
 \end{aligned}\right.
\end{equation}
where $[\trial]^{\nsd}:= \{\bu\in [\sobo(\Omega)]^{\nsd} : \bu=\bu_D \text{ on }\Gamma_D\}$. 
\SZ{Note that in the definition of $\bm{\trial}$ all the spaces in the tensorial product are raised to the power of $\nsd$, also the parametric ones. }
A standard weighted residuals approach, with integrals in $\Omega\times\bI$ and the usual integration by parts only in $\Omega$ produces a (spatially) weak form in this multi-dimensional setup. Namely, find $(\bu,p)\in\bm{\trial}\times\eltwo(\Omega\times\bI) $ such that
\begin{equation} \label{eq:Multi-D}
  A\bigl(\bu,\bv\bigr) + B\bigl(\bv,p\bigr) +B\bigl(\bu,q\bigr)
  =L\bigl(\bv\bigr),\;\forall (\bv, q)\in\bm{\mathcal{S}}\times\eltwo(\Omega\times\bI) ,
\end{equation}
where the test function space for velocities is $\bm{\mathcal{S}}:=[\test]^{\nsd}\otimes[\eltwo(\I_1)]^{\nsd}\otimes[\eltwo(\I_2)]^{\nsd}\otimes\dotsb\otimes[\eltwo(\I_{\npar})]^{\nsd}$ and
$[\test]^{\nsd}:= \{\bu\in[\sobo(\Omega)]^{\nsd} : \bu=\bm{0}\text{ on }\Gamma_D\}$.

The following definitions of the bilinear and linear forms are necessary:
\begin{equation} \label{eq:Forms2}
  \begin{split}
  A\bigl(\bu,\bv\bigr)&:= \int_{\I_1}\!\int_{\I_2}\dotsi\int_{\I_{\npar}} a\bigl(\bu,\bv\bigr) \, 
                                      d\mu_{\npar}\dotsm d\mu_2\, d\mu_1 ,\\
  B\bigl(\bu,q\bigr)&:= \int_{\I_1}\!\int_{\I_2}\dotsi\int_{\I_{\npar}} b\bigl(\bu,q\bigr) \, 
                                      d\mu_{\npar}\dotsm d\mu_2\, d\mu_1 ,\\
  L\bigl(\bv\bigr)&:=  \int_{\I_1}\!\int_{\I_2}\dotsi\int_{\I_{\npar}} 
                                 \ell\bigl(\bv\bigr) \,
                                 d\mu_{\npar}\dotsm d\mu_2\, d\mu_1, 
  \end{split}
\end{equation}
where
\begin{equation}\label{eq:Forms3}
 \begin{split}
 a\bigl(\bv,\bw\bigr) = \int_{\Omega}2\nu\,\bm{\nabla}\bv:\bm{\nabla}\bw\,d\Omega
 \, , &\quad
 b\bigl(\bv,q\bigr) = -\int_{\Omega}q\,\grad\cdot\bv\,d\Omega, \quad \textrm{ and}
 \\
\ell\bigl(\bv\bigr) = \int_{\Omega}\bm{s}\cdot\bv\,d\Omega 
                          + &\int_{\partial\Omega\setminus\Gamma_D} \bv\cdot\bm{t}\,d\Gamma.
 \end{split}
\end{equation}
	
Obviously, the number of dimensions of the solution domain increases with the number of parameters. To circumvent the \emph{curse of  dimensionality}, the PGD approach
\cite{Ammar-AMCR:06, Chinesta-CLBACGAAH:13, PGD-CCH:14, Chinesta-Keunings-Leygue}
is employed here. This approach assumes a separable structure in the approximation to $(\bu,p)$. Note that the tensor product spaces $\bm{\trial}$ and $\eltwo(\Omega\times\bI)$ inherit the multidimensional complexity of the problem and, in principle, do not assume separability of the functions. 	

\begin{remark}[Saddle point structure]\label{rm:saddle}
Note that equation \eqref{eq:Multi-D} is often written as 
\begin{equation} \label{eq:Saddle-Multi-D}
 \left\{\begin{aligned}
  A\bigl(\bu,\bv\bigr) + B\bigl(\bv,p\bigr) &= L\bigl(\bv\bigr)  &&\forall\bv\in\bm{\mathcal{S}} ,\\
                                   B\bigl(\bu,q\bigr) &= 0                      &&\forall q\in\eltwo(\Omega\times\bI) ,
 \end{aligned}\right.
\end{equation}
to evidence the saddle point problem at hand.
\end{remark}
\section{The Proper Generalized Approximation}\label{sec:PGD}
\subsection{Three alternative forms of the separated approximation}\label{sec:3alternatives}
As usual in a PGD strategy, a separated representation $(\upgd^n,\ppgd^n)$ is imposed to approximate the solution of \eqref{eq:Multi-D} in each subdomain. The couple $(\upgd^n,\ppgd^n)$ stands for the PGD approximation with $n$ terms (or modes) of the velocity-pressure couple and it is defined as a sum of separated terms. Each term (mode) is the product of functions depending only on one of the arguments $(\bx,\mu_1,\mu_2,\dotsc,\mu_{\npar})$. Note that, in some of the PGD implementations the separated modal functions are normalized and therefore a scalar coefficient affects each mode and characterizes its amplitude. The first \emph{mode}, $(\upgd^0,\ppgd^0)$, is arbitrarily chosen (for instance accounting for Dirichlet boundary conditions). Then, a greedy algorithm is implemented to compute successively the last one, that is to compute term $n$ assuming that term $n-1$ is available 
\cite{PGD-CCH:14, Chinesta-Keunings-Leygue}.

Three alternatives can be considered for this separation, see also \cite{Maday-MMPY:15}, depending on how the modal functions for parameters $\bmu$ are considered. 
\begin{description}
\item[Independent component-wise separation (case \#0):] 
A distinct parametric modal function is considered for each term of the PGD expansion providing $(\upgd^n,\ppgd^n)$,  each component of the velocity and the pressure ($\nsd+1$ components) and each parameter ($\npar$). Hence the total number of parametric modal functions is $n  (\nsd+1) \npar$. This is the most general separation because a different parameter function is considered for each velocity component and for pressure, namely, for $i=1, \dotsc,\nsd$
\begin{equation*}
 \left\{\begin{aligned}
 u_i(\bx,\mu) \approx u_{i_{\texttt{PGD}}}^n(\bx,\mu)
                       &= u_{i_{\texttt{PGD}}}^{n-1}(\bx,\mu) + F_{\bu,i}^{n}(\bx)\,
                             L_{\bu,i,1}^{n}(\mu_1)\, L_{\bu,i,2}^{n}(\mu_2)\dotsm L_{\bu,i,\npar}^{n}(\mu_{\npar}), \\
    p(\bx,\mu) \approx\ppgd^n(\bx,\mu)
                       &= \ppgd^{n-1}(\bx,\mu) + \Fp^{n}(\bx)\,
                             L_{p,1}^{n}(\mu_1)\, L_{p,2}^{n}(\mu_2)\dotsm L_{p,\npar}^{n}(\mu_{\npar})
 \end{aligned}\right.
\end{equation*}
\SZ{where $F$ denote the spatial functions (depending on $\bx$) characterizing the mode (subscripts indicate wether they refer to some component of the velocity or the pressure) and $L$ denote the parametric functions. All these are real-valued scalar fields. In the following the vector fields are boldfaced.}
This strategy imposes different parameter--dependent functions for each spatial velocity component. Consequently, the spatial differential operators are affected differently for each spatial component.  This induces enormous difficulties in order to express (at least well approximately) the different forms by the sum of products of parameter-dependent functions and parameter-independent operators. 

The major drawback of such an approach is that the incompressibility constrain will not be trivially enforced. That is, LBB or incompressibility stabilization must be specifically studied and it is not trivial due to the variable weightings introduced by the parameter functions. 

Moreover, the implementation of the alternated directions scheme in the PGD for this separation form is highly intrusive. This is because it requires a distinct treatment of the different directions, depending on the value of the assumed parametric modes. Thus, the use of non-intrusive strategies with commercial codes becomes much more involved.

\item[Unique parameter function for velocity independent]\textbf{from the ones for pressure (case \#1):}
With respect to the previous formulation, the functions affecting the $\nsd$ different components of the velocity are taken to be the same. Thus, the total number of parametric modal functions is equal to $n  \, 2 \, \npar$. Correspondingly, 
\begin{equation*}
 \left\{\begin{aligned}
 \bu(\bx,\bmu) \approx\upgd^n(\bx,\bmu)
                       &= \upgd^{n-1}(\bx,\bmu) + \Fu^{n}(\bx)\,
                             L_{\bu,1}^{n}(\mu_1)\, L_{\bu,2}^{n}(\mu_2)\dotsm L_{\bu,\npar}^{n}(\mu_{\npar}), \\
    p(\bx,\bmu) \approx\ppgd^n(\bx,\bmu)
                       &= \ppgd^{n-1}(\bx,\bmu) + \Fp^{n}(\bx)\,
                             L_{p,1}^{n}(\mu_1)\, L_{p,2}^{n}(\mu_2)\dotsm L_{p,\npar}^{n}(\mu_{\npar}). 
 \end{aligned}\right.
\end{equation*}
This approach uncouples naturally the parameter functions from the spatial ones when the spatial divergence is computed to impose incompressibility. This will have a major effect in simplifying the choice of the spatial spaces pairs for velocity and pressure. 
\begin{remark}[Divergence-free modes] \label{rm:DF}
%
Note that, given the separated representation of the PGD approximation, the point-wise divergence-free velocity, $\grad\cdot\upgd^n=0$, is guaranteed by imposing  $\grad\cdot\Fu^s=0$ for $s=0,\dotsc , n$. Likewise
\begin{equation*}\label{eq:weak-DF}
 b\bigl(\Fu^s,q\bigr) = 0  \quad\forall q\in\eltwo(\Omega)  \text{ and $s=0,\dotsc , n$}
 \Longrightarrow
 b\bigl(\upgd^n,q\bigr) = 0 \quad\forall q\in\eltwo(\Omega),
\end{equation*}
which imposes the usual weak divergence-free condition on $\Omega$ for any set of parameters $\bmu\in\bI$.
\end{remark}

\item[Unique parameter functions for velocity and pressure (case \#2):] 
Following with the simplification in the number of parameter functions at each mode, the next step is to employ the same function for every component of the velocity and also for pressure, with a total number of parametric modal functions equal to $n  \, \npar$, namely.
\begin{equation*}
 \left\{\begin{aligned}
 \bu(\bx,\bmu) \approx\upgd^n(\bx,\bmu)
                       &= \upgd^{n-1}(\bx,\bmu) + \Fu^{n}(\bx)\,
                             L_{1}^{n}(\mu_1)\, L_{2}^{n}(\mu_2)\dotsm L_{\npar}^{n}(\mu_{\npar}), \\
    p(\bx,\bmu) \approx\ppgd^n(\bx,\bmu)
                       &= \ppgd^{n-1}(\bx,\bmu) + \Fp^{n}(\bx)\,
                             L_{1}^{n}(\mu_1)\, L_{2}^{n}(\mu_2)\dotsm L_{\npar}^{n}(\mu_{\npar}). 
 \end{aligned}\right.
\end{equation*}
Obviously, this case also benefits from the separation of the divergence of the velocity as in the previous approximation.
\end{description}

As stated above, the first form of the separation (with different parameter modes affecting every component of the velocity and the pressure) leads to a cumbersome formulation requiring a highly intrusive implementation. Consequently, in the following, only the two latter alternatives are taken into consideration. 
For the sake of a simpler presentation and without any loss of generality, the subsequent developments are done for the particular case of only one parameter $\mu$ ($\npar=1$). Thus, the first alternative under consideration (case \#1) uses two parameter functions (one for $\bu$ and one for $p$), and the second alternative (case \#2) uses just one parameter function (the same for $\bu$ and $p$). In the general case of $\npar \ge 1$, the number of parameter functions in each case are $2\,\npar$ and $\npar$, respectively.

\subsection{Case \#1: Two parameter functions (one for $\bu$ and one for $p$)}\label{sec:case1}
For $\npar=1$, the PGD approximation is written in this case as
\begin{equation}\label{eq:sep}
 \left\{\begin{aligned}
 \bu(\bx,\mu) \approx\upgd^n(\bx,\mu)
                       &= \upgd^{n-1}(\bx,\mu) + \Fu(\bx)\, \Lu(\mu), \\
    p(\bx,\mu) \approx\ppgd^n(\bx,\mu)
                       &= \ppgd^{n-1}(\bx,\mu) + \Fp(\bx) \, \Lp(\mu).
 \end{aligned}\right.
\end{equation}
Note that, in order to shorten the writing, superscript $n$ is omitted in the notation of  the unknown functions $ \Fu^{n}(\bx)$, $\Lu^{n}(\mu)$, $\Fp^{n}(\bx)$ and  $\Lp^{n}(\mu)$.
\subsubsection{Solving for each mode}
The approximation defined in \eqref{eq:sep} is substituted in \eqref{eq:Multi-D} and tested in a tangent manifold. That is, the unknowns to be determined are $\Fu\in[\trial]^{\nsd}$, $\Lu\in\eltwo(\I)$, $\Fp\in\eltwo(\Omega)$, and $\Lp\in\eltwo(\I)$ such that 
\begin{equation} \label{eq:weak-mode}
    A\bigl(\Fu \, \Lu,\bv\bigr) + B\bigl(\bv,\Fp \, \Lp \bigr) + B\bigl(\Fu \, \Lu ,q\bigr)
 = R\bigl(\upgd^{n-1},\ppgd^{n-1},\bv\bigr) - B\bigl(\upgd^{n-1},q\bigr) ,
\end{equation}
for all $\bv$ and $q$ in the tangent manifold and 
being the residual $R(\cdot,\cdot,\cdot)$ defined by 
\begin{equation*}\label{eq:residual}
 R\bigl(\upgd^{n-1},\ppgd^{n-1},\bv\bigr)= L\bigl(\bv\bigr) - A\bigl(\upgd^{n-1},\bv\bigr) - B\bigl(\bv,\ppgd^{n-1}\bigr) .
\end{equation*}
The tangent manifold is readily characterized by choosing $\bv$ and $q$ as variations of $\Fu\,  \Lu$ and $\Fp \, \Lp$ respectively, that is
\begin{equation*}\label{eq:tangent}
 \bv=\de\Fu\,\Lu+\Fu\,\de\Lu \;\text{ and }\; q=\de\Fp\,\Lp+\Fp\,\de\Lp .
\end{equation*}
for all $\de\Fu\in[\test]^{\nsd}$, $\de\Lu\in\eltwo(\I)$, $\de\Fp\in\eltwo(\Omega)$, and $\de\Lp\in\eltwo(\I)$.

Following Remark \ref{rm:saddle}, this problem can also be equivalently rewritten as 
\begin{subequations}\label{eq:weak-mode-2}
  \begin{align}[left = {\empheqlbrace\,}]
A\bigl(\Fu \, \Lu ,\de\Fu\,\Lu\bigr) {+} B\bigl(\de\Fu\,\Lu,\Fp\, \Lp \bigr) 
  &= R\bigl(\upgd^{n-1},\ppgd^{n-1}, \de\Fu\,\Lu\bigr) &&  \forall \de\Fu\in[\test]^{\nsd} , \label{eq:weak-mode-2-1}\\
B\bigl(\Fu \, \Lu ,\de\Fp\,\Lp\bigr)  & = - B\bigl(\upgd^{n-1},\de\Fp\,\Lp\bigr) &&   \forall\de\Fp\in\eltwo(\Omega) ,  
\label{eq:weak-mode-2-2} \\
  A\bigl(\Fu \, \Lu ,\Fu\,\de\Lu\bigr) {+} B\bigl(\Fu\,\de\Lu,\Fp \, \Lp \bigr) 
  &= R\bigl(\upgd^{n-1},\ppgd^{n-1}, \Fu\,\de\Lu\bigr) &&\forall\de\Lu\in\eltwo(\I) , \label{eq:weak-mode-2-3} \\
  B\bigl(\Fu \, \Lu ,\Fp\,\de\Lp\bigr)  &= - B\bigl(\upgd^{n-1},\Fp\,\de\Lp\bigr)  &&\forall\de\Lp\in\eltwo(\I) . \label{eq:weak-mode-2-4}
\end{align}
\end{subequations}
Note that \eqref{eq:weak-mode-2}  is a nonlinear system of functional equations for the four unknowns $\Fu$, $\Lu$, $\Fp$, and $\Lp$. In the PGD framework, \eqref{eq:weak-mode-2}  is iteratively solved using an \emph{alternated directions} scheme. That is, first solving  \eqref{eq:weak-mode-2-1} and  \eqref{eq:weak-mode-2-2} for unknowns $\Fu$ and $\Fp$, assuming that $\Lu$ and $\Lp$ are known. This first stage is denoted \emph{spatial iteration} because it has the same structure of a standard (non parametric) Stokes problem. 

Then, equations \eqref{eq:weak-mode-2-3} and  \eqref{eq:weak-mode-2-4} are solved for unknowns $\Lu$ and $\Lp$ assuming that $\Fu$ and $\Fp$ are known. This step is denoted \emph{parameter iteration}, and it consists in iterating for every parametric direction (just once for $\npar=1$, in general $\npar$ steps are needed). The process is iterated between subsystem \eqref{eq:weak-mode-2-1} and  \eqref{eq:weak-mode-2-2} and 
subsystem \eqref{eq:weak-mode-2-3} and  \eqref{eq:weak-mode-2-4} until a stationary solution is reached.
 
 \begin{remark}[Solving groups of two equations] \label{rm:Two}
 In other PGD formulations, the alternated direction schemes for the nonlinear systems take the modes one by one, solving for one and assuming that the rest are known. Here, the two couples of unknowns $(\Fu, \Fp)$ and $(\Lu,\Lp)$ are solved together. This is due to the Saddle Point structure inherited by the groups of equations \eqref{eq:weak-mode-2-1} \&  \eqref{eq:weak-mode-2-2}  and \eqref{eq:weak-mode-2-3} \&  \eqref{eq:weak-mode-2-4}. In particular, the natural unknown for  \eqref{eq:weak-mode-2-2} would be $\Fp$ and it is not appearing explicitly in the equation. Thus, it is not possible solving  \eqref{eq:weak-mode-2-2} to find $\Fp$ assuming that $\Fu$, $\Lu$ and $\Lp$ are known. The same happens with $\Lp$ in  \eqref{eq:weak-mode-2-4}. 
  \end{remark}
\subsubsection{The spatial iteration.}\label{sec:SpaIter}
As stated above, the spatial iteration consists in solving  \eqref{eq:weak-mode-2-1} and  \eqref{eq:weak-mode-2-2} for unknowns $\Fu$ and $\Fp$, assuming that $\Lu$ and $\Lp$ are known.

The simplest separable form of the bilinear operators introduced in \eqref{eq:Forms2} is, for $\npar=1$
\begin{equation} \label{eq:Forms4}
  \begin{split}
  A\bigl(\Fu \Lu,\de\Fu \Lu \bigr)&= \int_{\I}\! \Lu \Lu a\bigl(\Fu,\de\Fu\bigr) \, d\mu  
                                      =  \left[ \int_{\I}  \Lu^{2}\, d\mu \right] a\bigl(\Fu,\de\Fu\bigr)  ,\\
  B\bigl(\de\Fu\,\Lu,\Fp\, \Lp)&= \left[ \int_{\I}  \Lu \, \Lp \, d\mu  \right] \, b\bigl(\de\Fu\,\Fp )  
  \end{split}
\end{equation}
In general, the separation of the bilinear form may require a sum of different terms. For the sake of a simple notation, this one-term separation is assumed to hold. The general case does not introduce additional conceptual complexity.

Thus, introducing the computable scalar quantities
\begin{equation}
 \alpha_{\mu}= \int_{\I} \Lu^2\,d\mu \; ,\qquad \beta_{\mu} = \int_{\I} \Lu\Lp\,d\mu \; , 
\end{equation}
the system of equations  \eqref{eq:weak-mode-2-1} and  \eqref{eq:weak-mode-2-2} reads
\begin{subequations}\label{eq:weak-mode-3}
  \begin{align}[left = {\empheqlbrace\,}]
 \alpha_{\mu} a\bigl(\Fu,\de\Fu\bigr)  + \beta_{\mu} b\bigl(\de\Fu\,\Fp )  
  &= R\bigl(\upgd^{n-1},\ppgd^{n-1}, \de\Fu\,\Lu\bigr) &&  \forall \de\Fu\in[\test]^{\nsd} , \label{eq:weak-mode-3-1}\\
  &=: \mathcal{R}_{\bu}\bigl(\upgd^{n-1},\ppgd^{n-1}, \de\Fu\,\Lu\bigr) && \nonumber \\
\beta_{\mu} b\bigl(\Fu,\de\Fp\bigr)  
& = - B\bigl(\upgd^{n-1},\de\Fp\,\Lp\bigr) &&   \forall\de\Fp\in\eltwo(\Omega)\label{eq:weak-mode-3-2}\\
  &=: \mathcal{R}_{p}\bigl(\upgd^{n-1},\de\Fp\,\Lp\bigr)   && \nonumber 
\end{align}
\end{subequations}
where the residual character of the left-hand-sides of \eqref{eq:weak-mode-3-1} and \eqref{eq:weak-mode-3-2} is emphasized introducing the notations $ \mathcal{R}_{\bu}$ and $\mathcal{R}_{p}$ such that for any $\bw\in [\test]^{\nsd}$ and $\omega \in\eltwo(\Omega)$, 
\begin{subequations}\label{eq:residuals}
\begin{gather}
  \begin{split}
  \mathcal{R}_{\bu}\bigl(\upgd^{n-1},\ppgd^{n-1}, \bw\,\omega\bigr) 
  = \biggr[\int_{\I}\omega \,d\mu\biggl] \, \ell\bigl(\bw\bigr) 
  - \sum_{s=0}^{n-1}&\biggr[\int_{\I}\omega\,\Lu^s \,d\mu\biggl]\,a\bigl(\Fu^{s},\bw\bigr)\\ 
  &- \sum_{s=0}^{n-1}\biggr[\int_{\I}\omega\,\Lp^s \,d\mu\biggl] \, b\bigl(\bw,\Fp^{s}\bigr), \label{eq:residuals1} 
  \end{split}\\
  \mathcal{R}_{p}\bigl(\upgd^{n-1},q\,\rho\bigr) 
  =   -\sum_{s=0}^{n-1}\biggr[\int_{\I} \Lu^s\,\rho\,d\mu\biggl] b\bigl(\Fu^s,q\bigr) . \label{eq:residuals2}
\end{gather}
\end{subequations}
Note that problem \eqref{eq:weak-mode-3} is linear for $\Fu$ and $\Fp$ and has the same structure of a standard (nonparametric) Stokes problem.

Once the discrete subspaces approximating $[\test]^{\nsd}$ and $\eltwo(\Omega)$ are chosen, the functional equation \eqref{eq:weak-mode-3} results in a linear system of algebraic equations. The matrix associated with the system in the Stokes model is symmetric, with $2\times 2$ blocks and a null submatrix on the diagonal, namely
\begin{equation*}
 \begin{pmatrix} \mat{K} & \mat{G}\\ \mat{G}^T & \mat{0} \end{pmatrix} .
\end{equation*}
A necessary condition to guarantee unicity of the solution is that the kernel of the gradient matrix $ \mat{G}$ reduces to the trivial space, that is $\ker\mat{G}=\{\vect{0}\}$, where $\ker\mat{G}:=\{\vect{q}: \vect{q}\in\RR^{\hat{n}} \text{ and } \mat{G}\vect{q}=\vect{0}\}$,  $\hat{n}$ being the number of pressure unknowns in the spatial domain. This implies that the standard finite element approaches for incompressibility can readily be applied in the context of the PGD parameterized Stokes problem. That is, the user-preferred choice of LBB spatial elements or incompressible stabilization can be directly used in this context.

\begin{remark}[Divergence-free solution]\label{rm:DF-sol}
Note that if the first term is weakly divergence free (for instance, this is trivial for homogeneous Dirichlet boundary conditions), namely $b\bigl(\Fu^0,\de\Fp\bigr)=0$ for all $\de\Fp\in\eltwo(\Omega)$, every mode will be weakly divergence free and consequently, following Remark \ref{rm:DF}, in this case $\upgd^n$ is weakly divergence-free. 
\end{remark}

%
\subsubsection{The parameter iteration}
Recall that this substep is made to determine the parameter functions for each mode and consists in solving \eqref{eq:weak-mode-2-3} \&  \eqref{eq:weak-mode-2-4} for $\Lu$ and $\Lp$ assuming that the spatial functions $\Fu$ and $\Fp$ are known. For the particular case of $\npar=1$, the problem is rewritten as 
\begin{subequations}\label{eq:mu-problem}
  \begin{align}[left = {\empheqlbrace\,}]
  \alpha_{\bu}\int_{\I}\de\Lu\, \Lu^{n}\, d\mu + \beta_{\bu} \int_{\I}\de\Lu\, \Lp^{n}\, d\mu
  &= \mathcal{R}_{\bu}\bigl(\upgd^{n-1},\ppgd^{n-1}, \Fu\,\de\Lu\bigr)   
  &&\forall\de\Lu\in\eltwo(\I) , \label{eq:mu-problem1} \\
  \beta_{\bu}\int_{\I}\de\Lp\, \Lu^{n}\, d\mu  &= \mathcal{R}_{p}\bigl(\upgd^{n-1}, \Fp\,\de\Lp\bigr)   
  &&\forall\de\Lp\in\eltwo(\I) ,\label{eq:mu-problem2} 
\end{align}
\end{subequations}
where
\begin{equation}\label{eq:defAuBu}
 \alpha_{\bu}= a\bigl(\Fu,\Fu\bigr) \; ,\qquad \beta_{\bu} =b\bigl(\Fu,\Fp\bigr)  \; ,
\end{equation}
and $\mathcal{R}_{\bu}$ and $\mathcal{R}_{p}$ defined in \eqref{eq:residuals} are the known separated expressions of the residuals for velocity and pressure at the previous PGD approximation $(\upgd^{n-1},\ppgd^{n-1})$ tested now with $\Fu\,\de\Lu$ and $\Fp\,\de\Lp$, respectively. 

\SZ{
\begin{remark}[Algebraic nature of  \eqref{eq:defAuBu}]\label{rm:AlgNat}
Note that equations \eqref{eq:defAuBu} for $\Lu$ and $\Lp$ are integral equations that do not derive from any differential equation but from algebraic ones. This can be readily shown by realizing that weighting function $\de\Lu$ and $\de\Lp$ could be taken (in a point collocation fashion) as a set of Dirac deltas, ensuring that the algebraic equation is fulfilled at all the points included in the collocation (the expressions in  \eqref{eq:defAuBu}  do not contain any derivative of the unknowns $\Lu$ and $\Lp$).
\end{remark}
}
\SZ{
In many PGD implementations, the parametric modes $\Lu$ and $\Lp$ are represented as Finite Element (FE) functional approximations (using the nodal values as degrees of freedom and the shape functions as functional basis) and equations \eqref{eq:defAuBu} are solved with a Galerkin approach (taking $\de\Lu$ and $\de\Lp$ equal to the shape functions), which in this case results in a standard Least Squares functional approximation. This is typically done in order to preserve in the computation of the parametric modes $\Lu$ and $\Lp$ the same coding structure as for the velocity and pressure modes  $\Fu$ and $\Fp$ when solving equation \eqref{eq:weak-mode-3}.

Thus,} once the discrete subspaces are chosen, a symmetric matrix is obtained, namely
\begin{equation*}
 \begin{pmatrix} \mat{M}_{\bu,\bu} & \mat{M}_{\bu,p}\\ \mat{M}_{\bu,p}^T & \mat{0} \end{pmatrix} .
\end{equation*}
Similarly as in the previous case, the condition that $\ker \mat{M}_{\bu,p}=\{\vect{0}\}$ ensures uniqueness of the solution. 

Note that, as stated in Remarks \ref{rm:DF} and \ref{rm:DF-sol}, the velocity modes are divergence-free ($\grad\cdot \Fu^{s}=0$ for $s=1,\ldots,n-1$) and therefore the right-hand-side of \eqref{eq:mu-problem2} is zero. This is provoking an inconsistency that is clearly demonstrated for the particular, but not at all unusual, case of using the same discrete subspace of $\eltwo(\I)$ for both $\Lu$ and $\Lp$. Under such an assumption, there is only one mass matrix $\mat{M} = \mat{M}_{\bu,\bu}=\mat{M}_{\bu,p}$, symmetric and positive definite (i.e. its kernel is zero). Thus problem \eqref{eq:mu-problem} has a unique solution and the system to solve for each parameter substep has the following structure:
\begin{equation*}
 \begin{pmatrix} \mat{M} & \mat{M}\\ \mat{M} & \mat{0} \end{pmatrix}
 \begin{pmatrix} \bm{L}_{u}   \\  \bm{L}_{p} \end{pmatrix} = 
  \begin{pmatrix} \bm{R}_{u}   \\  \bm{0} \end{pmatrix} .
\end{equation*}
This is obviously leading to an inconsistent solution of $\Lu=0$ and therefore $\upgd=0$.

This shows that the second approach for the separated representation is not viable.

\SZ{
Note that the same conclusion is reached by following a point collocation approach as described in Remark \ref{rm:AlgNat}.
Taking $\de\Lu=\de_{\mu}$  and $\de\Lp=\de_{\mu}$ in \eqref{eq:mu-problem2}, that is particularizing the algebraic equations 
for a given value of the parameter $\mu$, \eqref{eq:mu-problem2} results in
\begin{subequations}\label{eq:alg}
  \begin{align}[left = {\empheqlbrace\,}]
  \alpha_{\bu} \Lu^{n}(\mu) + \beta_{\bu}  \Lp^{n}(\mu)
  &= \mathcal{R}_{\bu}\bigl(\upgd^{n-1},\ppgd^{n-1}, \Fu\,\de_{\mu}\bigr)   
  && \label{eq:alg1} \\
  \beta_{\bu} \Lu^{n}(\mu)  &= \mathcal{R}_{p}\bigl(\upgd^{n-1}, \Fp\,\de_{\mu} \bigr)   
  && \label{eq:alg2} 
\end{align}
\end{subequations}
And, being  $\mathcal{R}_{p}\bigl(\upgd^{n-1}, \Fp\,\de_{\mu} \bigr)  = 0$  in \eqref{eq:alg2}, the solution is always $\Lu^{n}(\mu) = 0$. 

Thus, this drawback associated with the option taken in Case \#1 is independent of the choice of the functional description of the parametric modes and also of the approximation criterion to compute them.
}

\subsection{Case \#2: One parameter function (same for $\bu$ and $p$)}\label{sec:case2}
For $\npar=1$,  using a similar notation as in \eqref{eq:sep}, and with $L(\mu)$ replacing both $\Lu(\mu)$ and $\Lp(\mu)$ 
\begin{equation}\label{eq:sep2}
 \left\{\begin{aligned}
 \bu(\bx,\mu) \approx\upgd^n(\bx,\mu)
                       &= \upgd^{n-1}(\bx,\mu) + \Fu(\bx)\, L(\mu), \\
    p(\bx,\mu) \approx\ppgd^n(\bx,\mu)
                       &= \ppgd^{n-1}(\bx,\mu) + \Fp(\bx) \, L(\mu).
 \end{aligned}\right.
\end{equation}

The approximation defined in \eqref{eq:sep2} is substituted in \eqref{eq:Multi-D} and tested in a tangent manifold. Thus, the problem becomes,
find $\Fu\in[\trial]^{\nsd}$, $\Fp\in\eltwo(\Omega)$ and  $L\in\eltwo(\I)$ such that 
\begin{equation} \label{eq:weak-mode-final}
    A\bigl(\Fu \, L,\bv\bigr) + B\bigl(\bv,\Fp \, L \bigr) + B\bigl(\Fu \, L ,q\bigr)
 = R\bigl(\upgd^{n-1},\ppgd^{n-1},\bv\bigr) - B\bigl(\upgd^{n-1},q\bigr) ,
\end{equation}
for all $\bv$ and $q$ in the tangent manifold. Now, the space of unknowns and the tangent manifold have one dimension less with respect to the previous case, since $\Lu$ and $\Lp$ have been replaced by $L$. The corresponding expressions for the test functions are 
\begin{equation}\label{eq:tangent}
  \bv=\de\Fu\, L+\Fu\,\de L \;\text{ and }\; q=\de\Fp\, L+\Fp\,\de L .
\end{equation}

Thus, the equation corresponding to \eqref{eq:weak-mode-2} is derived by replacing also $\de \Lu$ and $\de \Lp$ by $\de L$, \begin{subequations}\label{eq:weak-mode-4}
  \begin{align}[left = {\empheqlbrace\,}]
A\bigl(\Fu \, L ,\de\Fu\, L \bigr) {+} B\bigl(\de\Fu\, L,\Fp\, L \bigr) 
  &= R\bigl(\upgd^{n-1},\ppgd^{n-1}, \de\Fu\, L\bigr) &&  \forall \de\Fu\in[\test]^{\nsd} , \label{eq:weak-mode-4-1}\\
B\bigl(\Fu \, L ,\de\Fp\, L \bigr)  & = - B\bigl(\upgd^{n-1},\de\Fp\, L \bigr) &&   \forall\de\Fp\in\eltwo(\Omega) ,  
\label{eq:weak-mode-4-2} \\
\end{align}
\begin{multline}
  A\bigl(\Fu \, L ,\Fu\,\de L \bigr) {+} B\bigl(\Fu\,\de L,\Fp \, L \bigr) {+} B\bigl(\Fu \, L ,\Fp\,\de L \bigr) \\
     = B\bigl(\Fu \, L ,\Fp\,\de L \bigr)  {-} B\bigl(\upgd^{n-1},\Fp\,\de L\bigr) \;\forall\de L\in\eltwo(\I) , \label{eq:weak-mode-4-3} 
\end{multline}
\end{subequations}
note that suppressing one unknown suppresses also one equation, because equations \eqref{eq:weak-mode-2-3} and  \eqref{eq:weak-mode-2-4} have been summed up into \eqref{eq:weak-mode-4-3}.

Again, equations \eqref{eq:weak-mode-4-1} and \eqref{eq:weak-mode-4-2} have to be solved together, due to the saddle point structure, as noted in Remark \ref{rm:Two}. 

The spatial iteration described in Section \ref{sec:SpaIter} is similar in this case. It consists in solving  \eqref{eq:weak-mode-4-1} and  \eqref{eq:weak-mode-4-2} for unknowns $\Fu$ and $\Fp$, assuming that $L$ is known. Recalling \eqref{eq:Forms4}
and introducing a new definition for the computable scalar quantity 
\begin{equation}
 \alpha_{\mu}= \int_{\I} L^2\,d\mu , 
\end{equation}
the system of equations  \eqref{eq:weak-mode-4-1} and  \eqref{eq:weak-mode-4-2} reads
\begin{subequations}\label{eq:weak-mode-5}
  \begin{align}[left = {\empheqlbrace\,}]
a\bigl(\Fu,\de\Fu\bigr)  +  b\bigl(\de\Fu\,\Fp )  
  &=  \mathcal{R}_{\bu}\bigl(\upgd^{n-1},\ppgd^{n-1}, \de\Fu\, L \bigr) /  \alpha_{\mu}&&  \forall \de\Fu\in[\test]^{\nsd} , \label{eq:weak-mode-5-1}\\
b\bigl(\Fu,\de\Fp\bigr)  
& = \mathcal{R}_{p}\bigl(\upgd^{n-1},\de\Fp\, L \bigr) /  \alpha_{\mu} &&   \forall\de\Fp\in\eltwo(\Omega)\label{eq:weak-mode-5-2}
\end{align}
\end{subequations}
Note that problem \eqref{eq:weak-mode-5} has the same structure as problem \eqref{eq:weak-mode-3} and therefore the spatial iterations are equivalent for cases \#1 and \#2.

The structure of the parameter iteration, however, changes considerably in case \#2. It consists in solving \eqref{eq:weak-mode-4-3} for $L$ assuming that $\Fu$ and $\Fp$ are known. Recalling the definitions in \eqref{eq:defAuBu}, \eqref{eq:weak-mode-4-3} becomes a unified version of the equation in \eqref{eq:mu-problem} (the sum of \eqref{eq:mu-problem1} and \eqref{eq:mu-problem2}), namely
\begin{equation}\label{eq:mu-problemS}
(\alpha_{\bu} + 2 \beta_{\bu} ) \int_{\I}\de L \, L \, d\mu 
= \mathcal{R}_{\bu}\bigl(\upgd^{n-1},\ppgd^{n-1}, \Fu\,\de L \bigr) 
+ \mathcal{R}_{p}\bigl(\upgd^{n-1}, \Fp\,\de L \bigr) \quad \forall\de L \in\eltwo(\I) .
\end{equation}
Note that this problem results in a simple system of equations with just a mass matrix, with the right-hand-side accounting for the effect of all the residuals. The system is easily solvable and provides a single parametric mode, affecting both velocity and pressure modes, with no particular restrictions.

Consequently, the alternative analyzed as case \#2 appears to be viable and, as confirmed in the numerical examples shown in the next section, is the right approach to define a parametric separation of the saddle point problems.

\subsection{Least-squares PGD projection and PGD compression.}\label{sec:PGDcompression}

Often, the PGD separated solution is post-processed with a compression algorithm based on a least squares projection in order to reduce the number of PGD modes. This is standard in the PGD practice, because the PGD terms may contain some intrinsic redundancy that is alleviated with this post-process. When compared with an SVD separation of the complete parametric solution, the redundancy is associated with the nonorthogonality of the different terms (or, conversely, the optimality of the SVD representation is associated with their orthogonality).

The projection strategy is described in \cite{DM-MZH:15} where its superior performance is demonstrated when compared with standard SVD (for a 2D separation) and HOSVD (High-Order SVD, for a larger number of parameters). It consists in finding with a PGD like algorithm (that is, greedy and with an alternated directions iterative scheme) the separated functions $\uupgd$ and $\pppgd$ that better approximate $\upgd$ and $\ppgd$ with a least squares criterion. Namely, $\uupgd$ and $\pppgd$ are sought such that 
\begin{equation} \label{eq:LSprojection}
    \bigl(\upgd^n - \uupgd,\bv\bigr) = 0 \text{ and }  \bigl(\ppgd^n-\pppgd ,q\bigr) = 0,
\end{equation}
for all $\bv$ and $q$ ranging in some suitable spaces. Equations \eqref{eq:LSprojection} are solved with a PGD strategy, exactly as described in the previous sections for the Stokes problem and equation \eqref{eq:Multi-D}. Typically, this operation is performed selecting $\uupgd$ and $\pppgd$ in the same functional  spaces as $\upgd$ and $\ppgd$. Bur here, we consider also projecting into richer functional spaces.

In the present context, this technology deserves a particular attention because it will allow addressing a concern that is naturally raised after the considerations introduced in the previous sections. A clear conclusion of the above analysis is that the PGD solution of the Stokes problem must adopt the formulation labelled as case \#2 (same parametric mode for all the velocity components and the pressure), while case \#1 (same parametric mode for all the velocity components and a different one for the pressure) is not viable. Moreover, the first idea announced  in section \ref{sec:3alternatives} that we may label now as case \#0 (all parametric modes different) was also discarded because of its implementation complexity in commercial codes.

The discarded forms (cases \#0 and \#1) are richer descriptions of the solution, in the sense that, for the same number of terms in the PGD sum, the number of degrees of freedom used to describe $\upgd$ and $\ppgd$ is much larger in cases \#0 and \#1 than in case \#2. Roughly speaking, the number of d.o.f. describing the parametric modes is in case \#1 the double of case \#2 and the factor goes to $d+1$ (being $d$ the number of spatial dimensions) when compared to case \#0.  Considering the total amount of d.o.f., one may think that the \emph{richer} descriptions would require less terms in the PGD sum, with the subsequent computational savings. Thus, the question is: does the proposed PGD formulation of the Stokes problem requires an excessive number of PGD terms?

The PGD least squares projection is used here to answer this question. Indeed, the three alternatives are admissible and viable to solve equation \eqref{eq:LSprojection} with a PGD approach. 
Thus, once $\upgd$ and $\ppgd$ are computed as described in case \#2,  $\uupgd$ and $\pppgd$ may  be computed solving equation \eqref{eq:LSprojection} with any of the three formulations. This will allow checking if having more d.o.f. per term in the sum results in having a shorter PGD sum (less terms). In other words, it will indicate if the restriction of having the same parametric mode for all the components of velocity and the pressure is an artificial constraint or if, on the contrary, this additional condition fits the form of the parametric solution. The numerical examples presented in the next section demonstrate that, in a pretty general situation, the case \#2 formulation does not significantly increase the number of terms in the PGD sum with respect to the alternatives corresponding to cases \#0 and \#1.

\SZ{
In short, case \#1 uses $\Lu \neq \Lp$ and case \#2 $\Lu = \Lp$ and consequently one would expect that one term of case \#1
\begin{equation*}
 \begin{pmatrix} \Fu \Lu \\ \Fp \Lp \end{pmatrix} =   \begin{pmatrix} \Fu^{1} L^{1} \\ \Fp^{1} L^{1} \end{pmatrix} + 
 \begin{pmatrix} \Fu^{2} L^{2} \\ \Fp^{2} L^{2} \end{pmatrix}
\end{equation*}
requires two terms of case \#2 (for instance, taking $ \Fu^{1}=\Fu$,  $\Fp^{1}=0$, $L^{1}=\Lu$, $ \Fu^{2}=0$,  $\Fp^{2}=\Fp$ and $L^{2}=\Lp$).

In other words, we claim that the solutions of the examples analyzed in the next section are optimally represented by the case \#2  option, in the sense that reducing the number of degrees of freedom per PGD term does not increase the number of PGD terms required. 

Since the alternative solutions (with cases \#0 and \#1) of the original problem are discarded, this assertion is demonstrated by representing (using the PGD compression) the solution obtained with case \#2 in the forms of cases \#0 and \#1. Instead of saving a significant number of PGD terms, the compression obtained is not relevant: the number of terms required is very similar. Thus, we conclude the solution fits naturally with the functional structure provided in case \#2 and the fact of reducing the number of degrees of freedom per PGD term it is not introducing extra terms in the PGD solution. 
}

\section{Numerical examples} \label{sec:exa} 

This section presents three examples, the first two are a parametric Stokes problems, with  a set of parameters determining the geometry of the computational domain. The first example is a very simple backward facing step with a single geometric parameter. This example is used to discuss the different alternatives for the PGD representation (cases \#1 and \#2 described in sections \ref{sec:case1} and \ref{sec:case2}) and to analyze the effect of the corresponding PGD compression techniques. 

The second example considered in this section describes a Stokes flow around a NACA airfoil, where the geometry depends on four independent parameters. 
The analysis includes a discussion on the how the PGD compression techniques are affected by the choice of the hybrid formulation. First, the particular treatment of geometric parameters in the PGD formulation is briefly recalled.

The third example addresses the Brinkman problem with a free parameter stating the relative weight of the Stokes and Darcy models (Brinkman is seen as a combination of both). As a consequence of the conclusion of the first example and the previous section, among the PGD formulations discussed above, only the alternative \#2 (see section \ref{sec:case2}) is considered in the second example.

\subsection{Stokes flow in domains with parametric geometry}

\subsubsection{Accounting for geometric parameters}
The strategy to deal with geometric parameters in the PGD solver was devised in  \cite{AH-AHCCL:14} for Poisson problems and then combined with material parameters in \cite{SZ-ZDMH:15} and \cite{NME:NME5313} for heat and wave propagation problems. The fundamental idea is using a reference domain $\mathcal{T}$ and a parametric mapping to the physical domain $\Omega({\bm \mu})$. Thus, the (physical) problem is stated in the reference domain. The problem in the reference domain includes some fictitious parametric properties accounting for the mapping. The mapping between $\mathcal{T}$ and $\Omega({\bm \mu})$ is described with a coarse FE mesh, much coarser than the computational mesh because it is only required to resolve the parametric variations of the geometry. 

In the case of a Stokes problem \eqref{eq:Stokes}, the bilinear form in the left-hand-side of the weak equation reads
\begin{align*}
 a(\bu,\bv) &=\int_{\Omega({\bm \mu}) } \grad \bu \cdot (k \grad \bv) \,d\Omega \\
           &= \int_{\mathcal{T}} \grad_{\!\hat{\bm x}} \bu \cdot 
                (\underbrace{ k \, | {{\bf J}({\bm \mu})} | {{\bf J}({\bm \mu})}^{-\textsf{T}}  {{\bf J}({\bm \mu})}^{-1} }_{{\bf D}({\bm \mu})}   \grad_{\!\hat{\bm x}} \bv)  \,d\hat{\bm x}
\end{align*}
where ${\bf J}({\bm \mu})$  is the Jacobian of the mapping and $\hat{\bm x}$ are coordinates in the reference domain  $\mathcal{T}$. The main goal of using this mapping is to transform the parametric dependence of the integration domain $\Omega({\bm \mu})$ into a parametric dependence of material-like properties, the parametric fictitious conductivity ${\bf D}({\bm \mu})$. In order to use PGD, the parametric dependence of ${\bf D}({\bm \mu})$ has to be expressed (often it has to be approximated) in a separable form. This requires a further step computed via SVD or higher-order SVD as described in detail in \cite{SZ-ZDMH:15}.

\subsubsection{Backward facing step}
\begin{figure}
 \centering
 \includegraphics[width=0.75\textwidth]{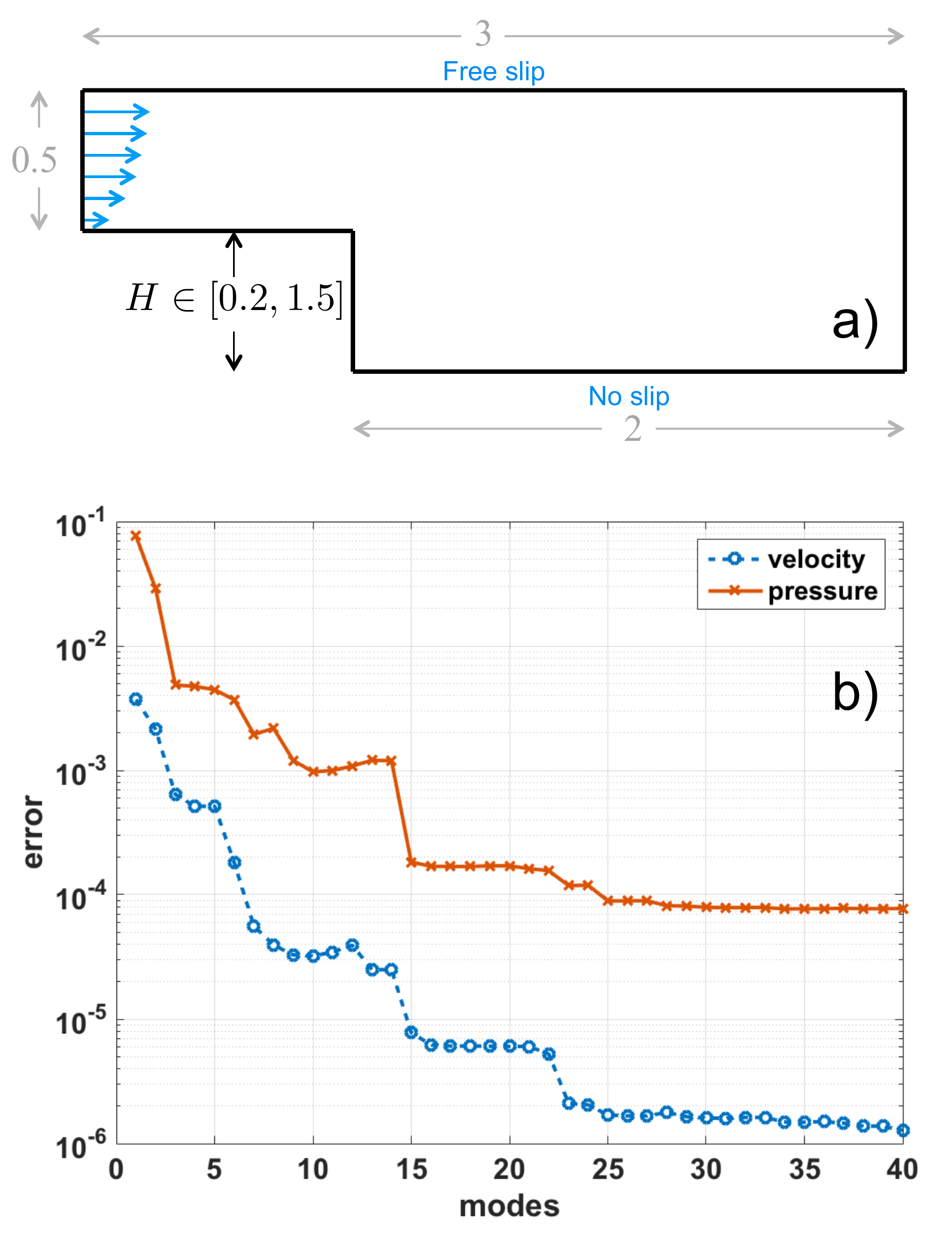}
 \caption{ Backward facing step example. 
 (a) \SZ{Setup of the model: domain and channel size is indicated in the Figure. Velocity boundary conditions are free slip on the top wall, a parabolic velocity profile on the inlet, no slip in the bottom wall (including the bottom part of the channel and the vertical wall of the step and Neumann homogeneous in the outlet. No pressure boundary conditions are required.}
 (b) Evolution of the error of the PGD solution (velocity and pressure) with the number of modes (maximum error compared with the FE solution for all the parametric values in a grid discretizing the parametric space).}
 \label{fig:bfs}
\end{figure}
\begin{figure}
 \centering
 \includegraphics[width=1\textwidth]{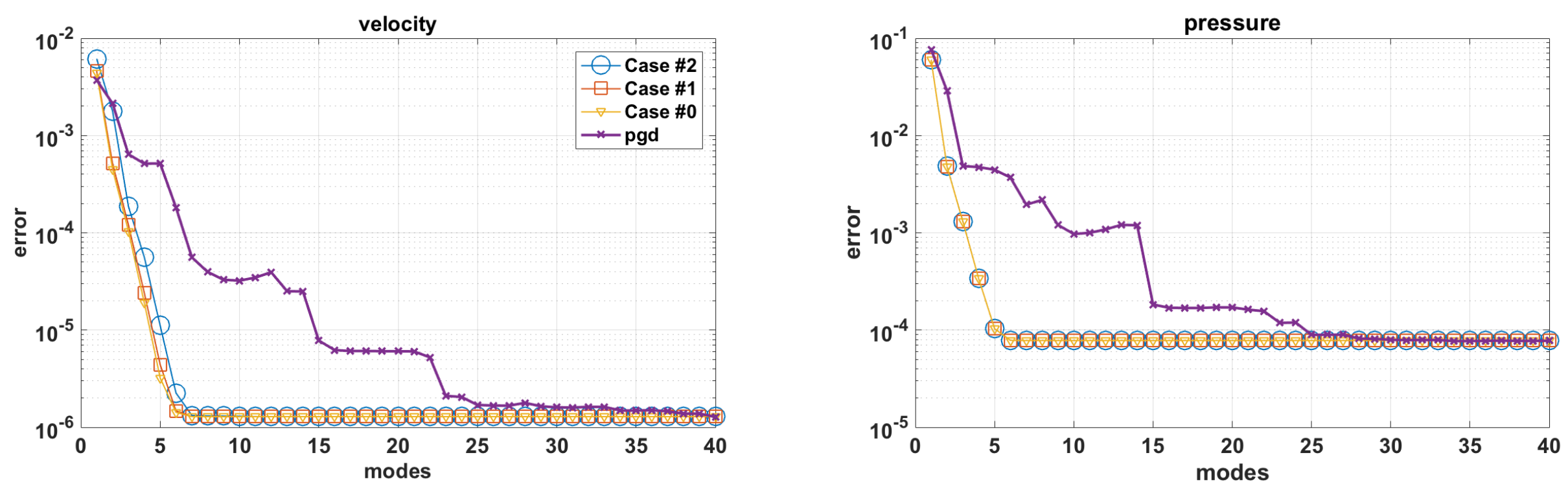}
 \caption{\SZ{Backward facing step example. Evolution of relative error \SZ{(measured in $\eltwo$ norm)} with the number of modes of the velocity and pressure fields corresponding to the PGD solution and the three Least-Squares projections (compressions) using the different formulations.}}
 \label{fig:projection:error}
\end{figure}

\SZ{A simplified backward facing step problem based on the Stokes equation}, see Figure~\ref{fig:bfs}, is considered where the parameter $H$ is the height of the step. The evolution of the error of the PGD  in terms of the number of modes is also shown in Figure~\ref{fig:bfs}. The error is taken as compared with a Finite Element solution for specific values of the parameter, measured in the supreme norm, that is the infinite norm, for the space dimension. Note that the error decreases with the number of terms and is larger for the pressure than for the velocity where it stagnates at a relative error of $10^{-6}$ for 40 PGD terms.

The same curve is repeated in Figure~\ref{fig:projection:error}, this time accompanied by three others representing the errors of the compressed PGD solutions. The additional three curves correspond to Least Squares projections using the three different alternatives, as described in Section~\ref{sec:PGDcompression}. Recall that the three alternatives correspond to cases \#0, \#1 and \#2 with different parametric modes in the solution. The case \#0 defines the larger functional space (number of parametric modes equal to number of space dimensions plus one) which includes the functional space of  case \#1 (number of parametric modes equal to two). Case  \#2 defines a smaller functional space with a single parametric mode. In consequence, the lower projection error must correspond to case \#0, then the error of case \#1 is larger (or equal) and even larger (or equal) for case \#2.

The results demonstrate that the three projection behave almost equally. The inclusion of the functional spaces is indeed translated into the expected inequality of the errors (error \#0 $\le$ error \#1 $\le$ error \#2 ) but the difference is very small. This reveals that selecting option \#2 is not practically increasing the required number of PGD terms with respect to a richer functional description (being case \#0 the richest). The structure of the parametric description enforced in case \#2 seems to fit the nature of the actual solution and its parametric dependence.

\subsubsection{Stokes flow around a NACA airfoil}

\begin{figure}
 \centering
 \includegraphics[width=0.75\textwidth]{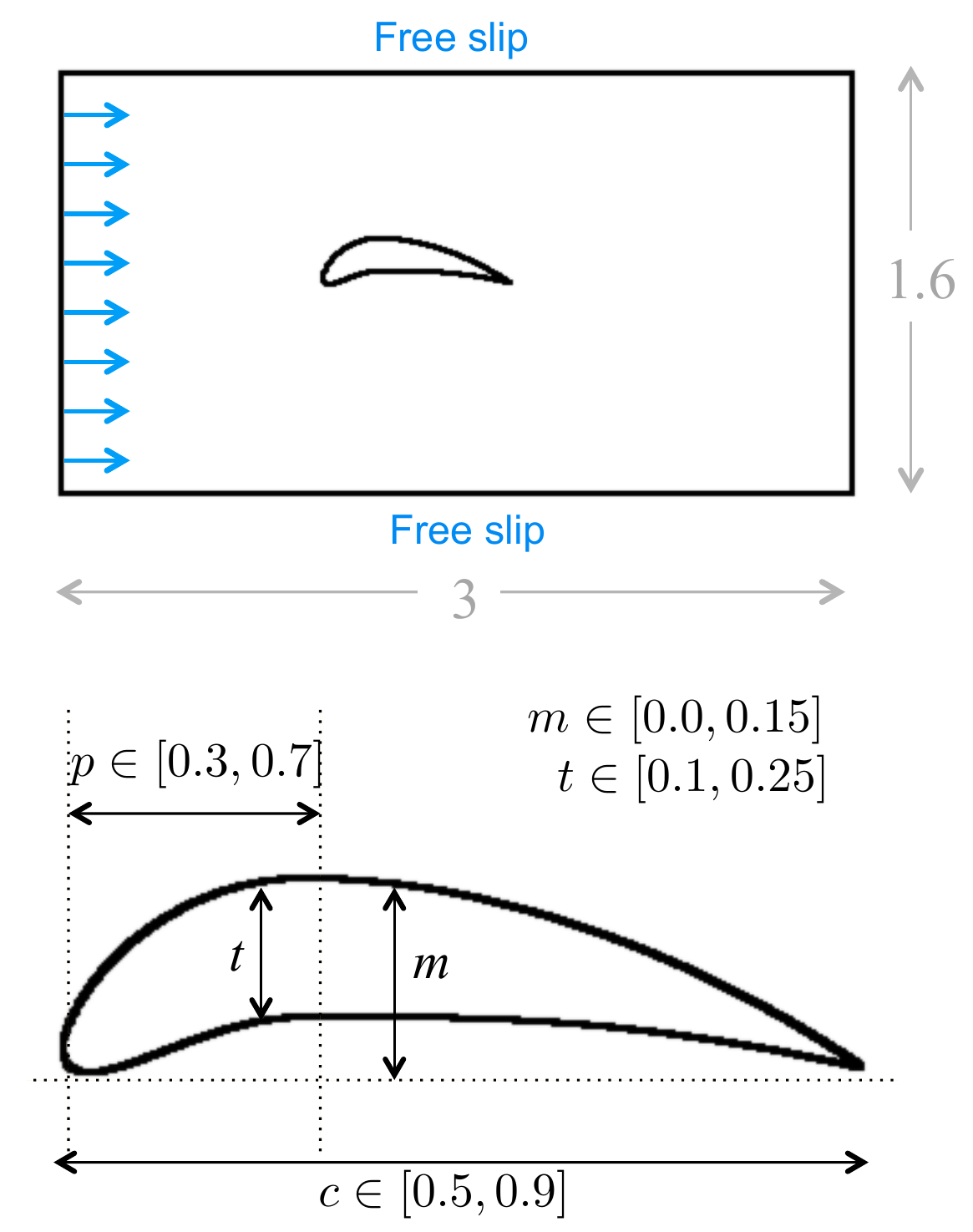}
 \caption{\SZ{Stokes flow around a 4-digits NACA airfoil parameterized on 4 quantities describing the airfoil geometry. 
Velocity  ${\bm u}(x,c, t, m, p)$ depends on space and the following parameters:
the chord length $c\in[0.5, 0.9]$,
the maximum thickness as a fraction of the chord $t\in[0.1, 0.25]$,
the maximum camber $m\in[0, 0.15]$, and 
the location of maximum camber $p\in[0.3, 0.7]$. All parameters are discretized using 20 linear elements.
The size of the computational domain is described in the Figure. Velocity boundary conditions are free slip on horizontal walls, a constant horizontal velocity with value one is imposed on the inflow wall and Neumann homogeneous conditions on the outflow. No pressure boundary conditions are required.}}
 \label{fig:naca:param}
\end{figure}

The flow around a 4-digit NACA airfoil is considered, see Figure \ref{fig:naca:param}. The four geometric parameters (digits, in this contex) defining the geometry of the airfoil are: 1) length, $c$, 2) thickness, $t$, 3) max camber, $m$ and  4) max camber position, $p$. 

\begin{figure}
 \centering
 \includegraphics[width=0.7\textwidth]{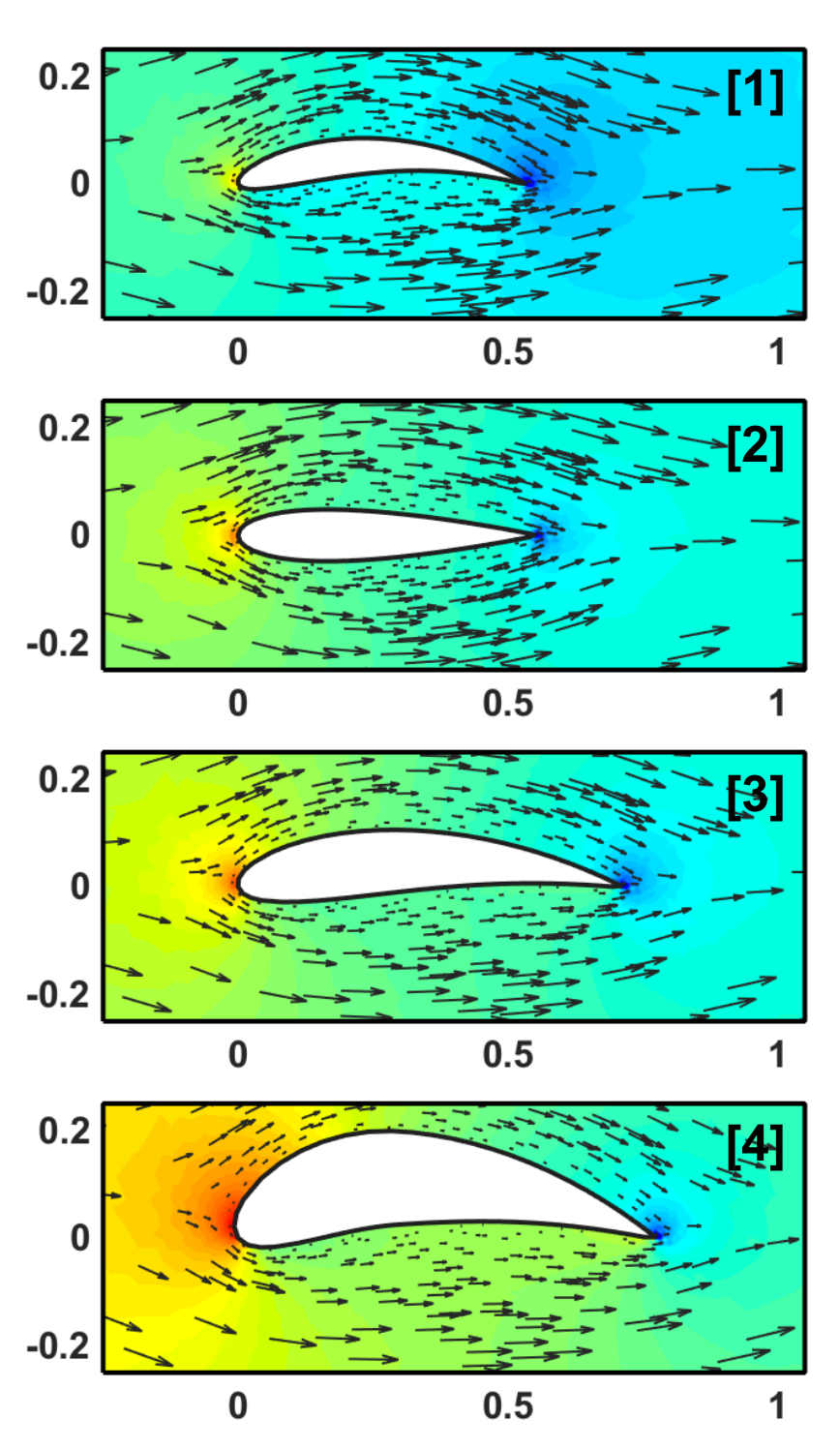}
 \caption{NACA-4 example. Velocity (quiver fields) and pressure (color maps) solutions of the parametrized NACA-4 airfoils for four representative sets of parameter values\SZ{  \{c,t,m,p\}:
 [1]=\{0.53, 0.13, 0.10, 0.51\};
 [2]=\{0.55, 0.17, 0.00, 0.50\};
 [3]=\{0.71, 0.17, 0.07, 0.51\};
 [4]=\{0.77, 0.23, 0.14, 0.38\}.
 }} 
 \label{fig:naca:res}
\end{figure}
\begin{figure}
 \centering
 \includegraphics[width=0.75\textwidth]{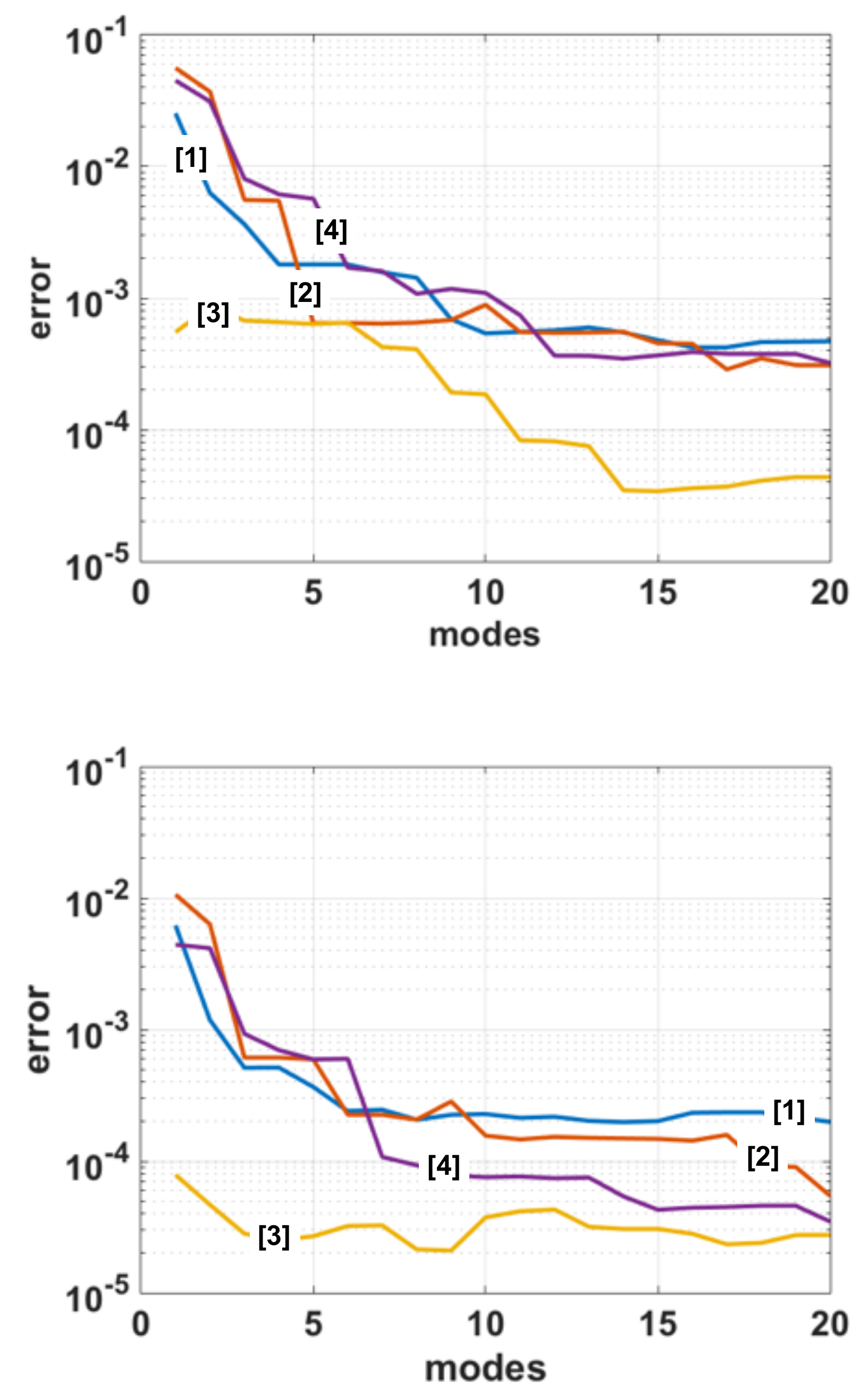}
 \caption{NACA-4 example. Evolution with the number of PGD modes (in abscissae) of the relative error for velocity (top, measured in $H^{1}$ norm) and pressure (bottom, measured in $L_{2}$ norm) for the particular parametric values selected in Figure \ref{fig:naca:res}. The PGD error is computed with respect to the complete FE solution obtained for the particular parametric values.}
 \label{fig:naca:conv}
\end{figure}

The PGD solution provides a computational vademecum containing the flow solutions for any possible NACA-4 geometry. For the sake of illustration,  the evaluation for four particular values of the parameters is shown in Figure~\ref{fig:naca:res}. These four particular parametric values are also used to check convergence with the number of PGD modes. Figure~\ref{fig:naca:conv} presents the evolution of the errors in velocity and pressure as the number of modes increases. 

Note that also in this complex example (with four parametric dimensions) the PGD solution behaves correctly when compared to standard FE solutions for specific (and representative) values of the parameters. None of this four sampling points of the 4D parametric space is a grid point (the values of the parameters selected do not coincide with the discrete grid of each parametric dimension). Note that integration in the parametric space to compute a $L_{2}$ norm of the error is avoided because of the associated computational burden.

\subsection{Flow in fractured media}

\begin{figure}
 \centering
 \includegraphics[width=0.75\textwidth]{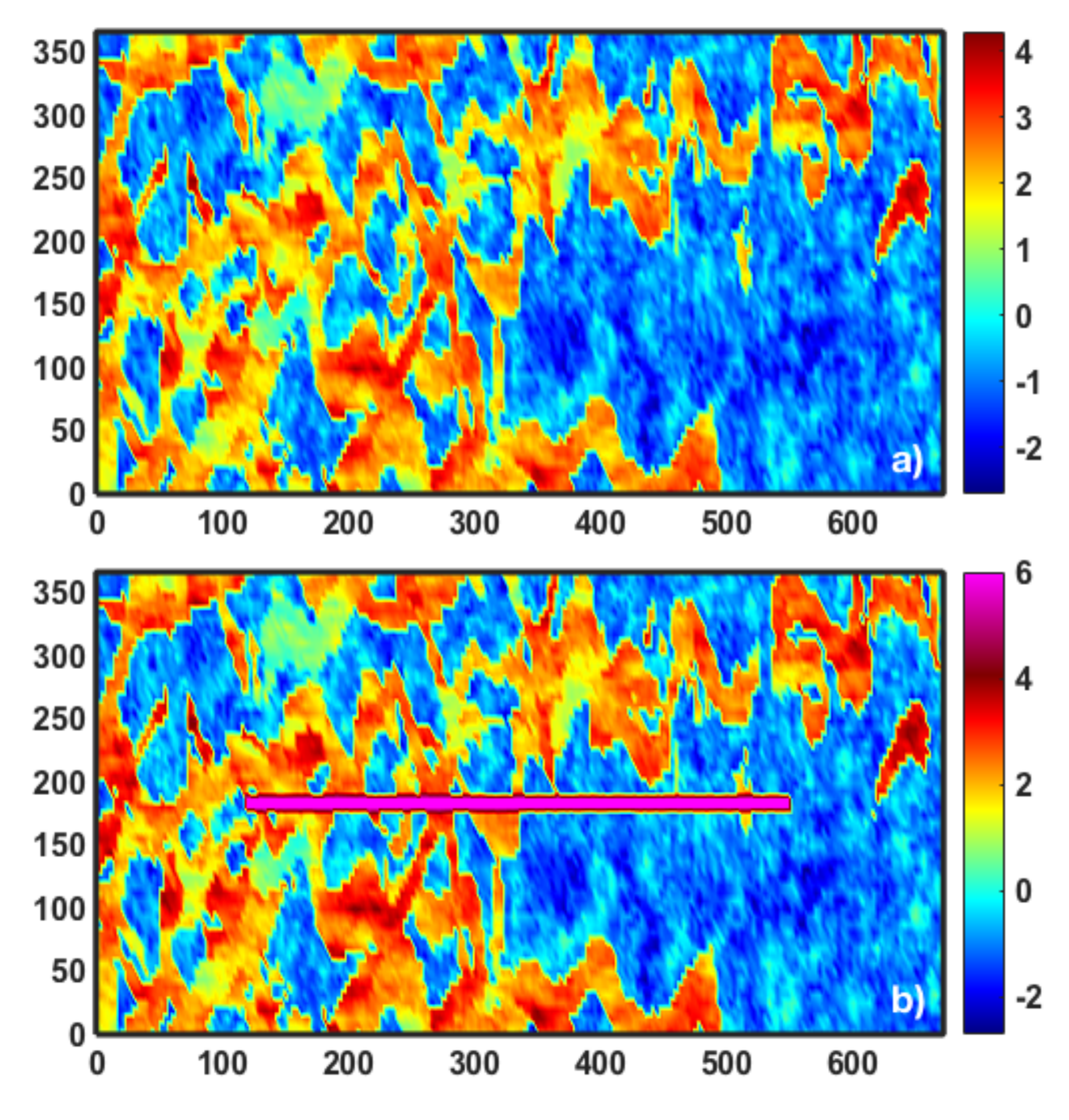}
 \caption{Fractured medium. Spatial distribution of the original permeability $k_{\text{SPE10}}(\bx)$ \cite{Christie:01} (top) and the perturbed parametric permeability accounting for the fracture $k(\bx,\mu_{\text{max}}=6)$ (bottom). The scales are logarithmic, ($\log_{10}$)  permeabilities in milidarcy and lengths in meters along axes. 
\SZ{The computational domain is 671$\times$336 m in size and it is discretized in 220$\times$60 quadrilateral elements with order 2 for velocities and order 1 for pressures. The size of the horizontal channel is 430$\times$11.2 m and its lower left corner is located in (122, 168) m. The flow is driven by an imposed pressure of 1 cP on the point (0,213) m.}}
 \label{fig:permeability}
\end{figure}

\begin{figure}
 \centering
 \includegraphics[width=0.65\textwidth]{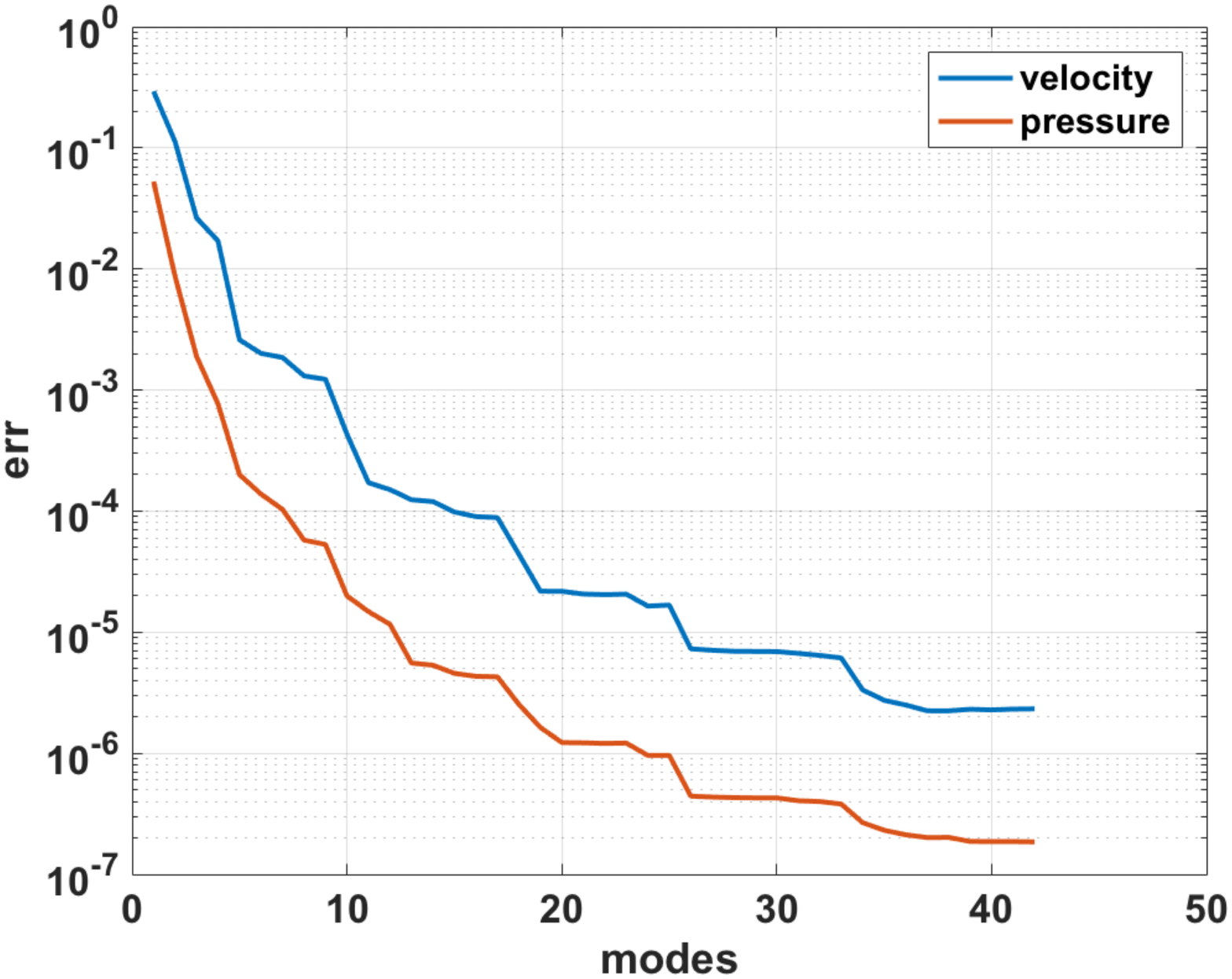}
 \caption{Fractured medium. Error of the PGD solution as a function of the number of PGD terms. The relative error is computed 
using the maximum of the difference between PGD and FE solutions (infinite-norm; FE solutions are calculated  for every parameter value).}
 \label{fig:brinkman:convergence}
\end{figure}

Describing flow in fractured, vuggy and porous media with a unique model is important for reservoir engineering.
A standard approach is using Darcy's law in the porous domain and Stokes law in the parts of the domain containing voids and fractures. 
The Brinkman model merges both Darcy and Stokes, see \cite{Brinkman:47,Popov:09,Konno:12}, linearly combining the effect of the two constitutive models. Thus, the Brinkman problem reads: find  velocity $\bu$ and pore pressure $p$ such that
\begin{subequations}\label{eq:brinkman}
  \begin{align}[left = {\empheqlbrace\,}]
-\tilde \eta \Delta \bu  + \eta \mat{K}^{-1} \bu + \nabla p &= 0&&\text{in $\Omega$} \label{eq:brinkman1} \\
\nabla \cdot \bu &= 0&&\text{in $\Omega$}, \label{eq:brinkman2}
\end{align}
\end{subequations}
where $\eta$ and $\tilde \eta$ are the dynamic and effective viscosity respectively and $\mat{K}$ is the permeability tensor. Note that the dynamic viscosity $\eta$ allows tuning the weight of the Darcy model in the constitutive description of the flow. For small values of $\eta$, a free-flow Stokes-like pattern is obtained and  for large values of $\eta$ the solution tends to behave as porous flow (Darcy). In reservoir modelling it is usual to assume incompressibility \eqref{eq:brinkman2}, to neglect gravity (right-hand-side of \eqref{eq:brinkman1}) and to set $\eta=\tilde \eta$ for the bulk material, \cite{Popov:09,Konno:12}.

A parametrized version of the Brinkman problem is built upon the background permeability field provided as a test case (SPE10) by the Society of Petroleum Engineers  \cite{Christie:01}. The original three-dimensional (3D) field is restricted to a 2D domain following K\"onn\"o and Stenberg \cite{Konno:12}. Here, we aim at analyzing the effect of adding a fracture to the layer 68 of the SPE10 model. Note that SPE10 provides a non-uniform isotropic permeability and, therefore, the matrix $\mat{K}$ becomes a scalar field $k_{\text{SPE10}}(\bx)$. In the framework of the Brinkman model, the fracture is accounted for in a natural way by significantly increasing the permeability in the region where the fracture is located. 
Figure \ref{fig:permeability} shows the SPE10 permeability with and without the perturbation that accounts for the fracture. 

The value of the perturbed permeability in the fractured zone is not easy to set. The modeller knows only that permeability has to be significantly larger in the zones where the free-flow pattern is expected, with respect to the reference values for porous flow. Thus, it is particularly interesting having a tuning parameter that allows enforcing a gradual transition between the two regimes. This is especially relevant in flow simulations of karst reservoirs where vugs and caves are embedded in a porous rock and are connected via fracture networks at multiple scales. For example the work of Popov and coauthors \cite{Popov:09} based on a Brinkman model, assumes a continuous permeability ranging six orders of magnitude and investigates the effect of permeability in a filled fracture.

\begin{figure*}
 \centering
 \includegraphics[width=\textwidth]{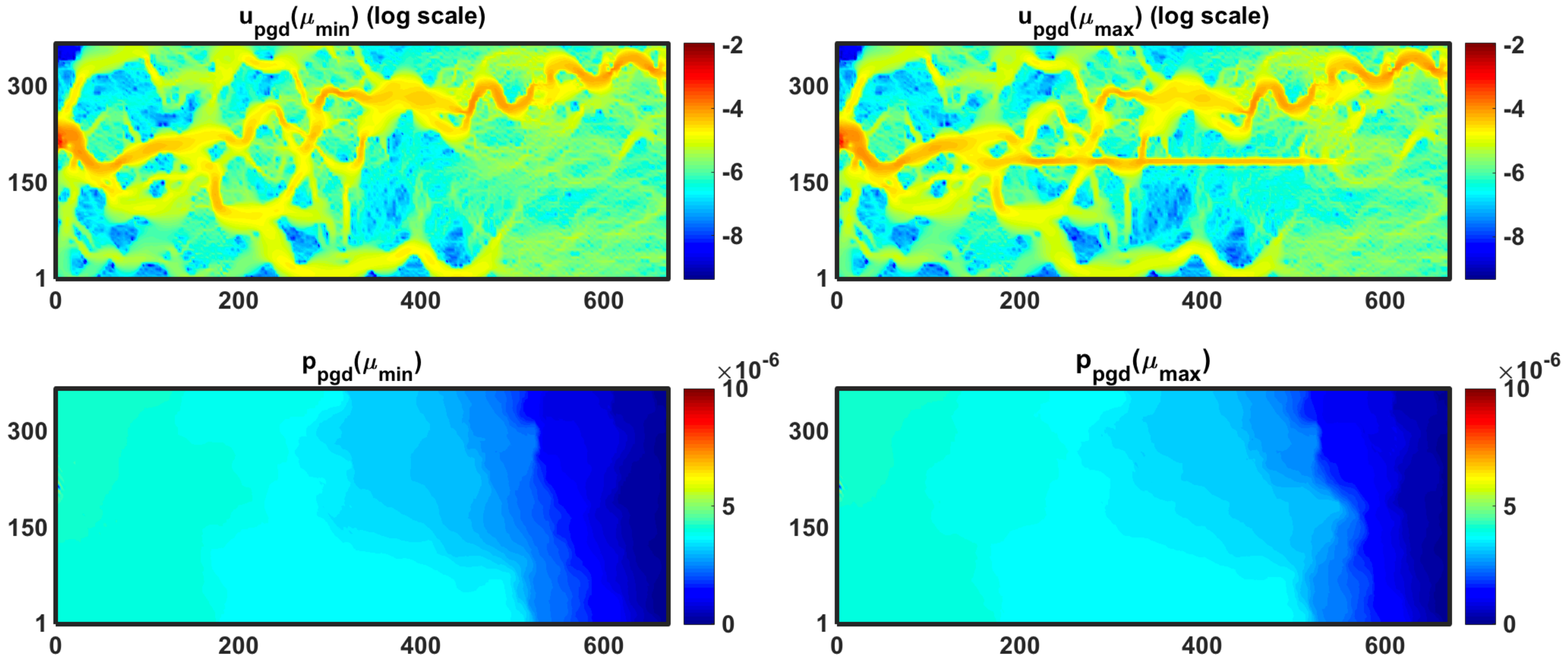}
 \caption{Norm of the velocity and pressure provided by PGD  (top and bottom rows respectively) evaluated at minimum and maximum value for the parameter (left and right columns respectively).}
 \label{fig:brinkman:solution}
\end{figure*}

\begin{figure*}
 \centering
 \includegraphics[width=\textwidth]{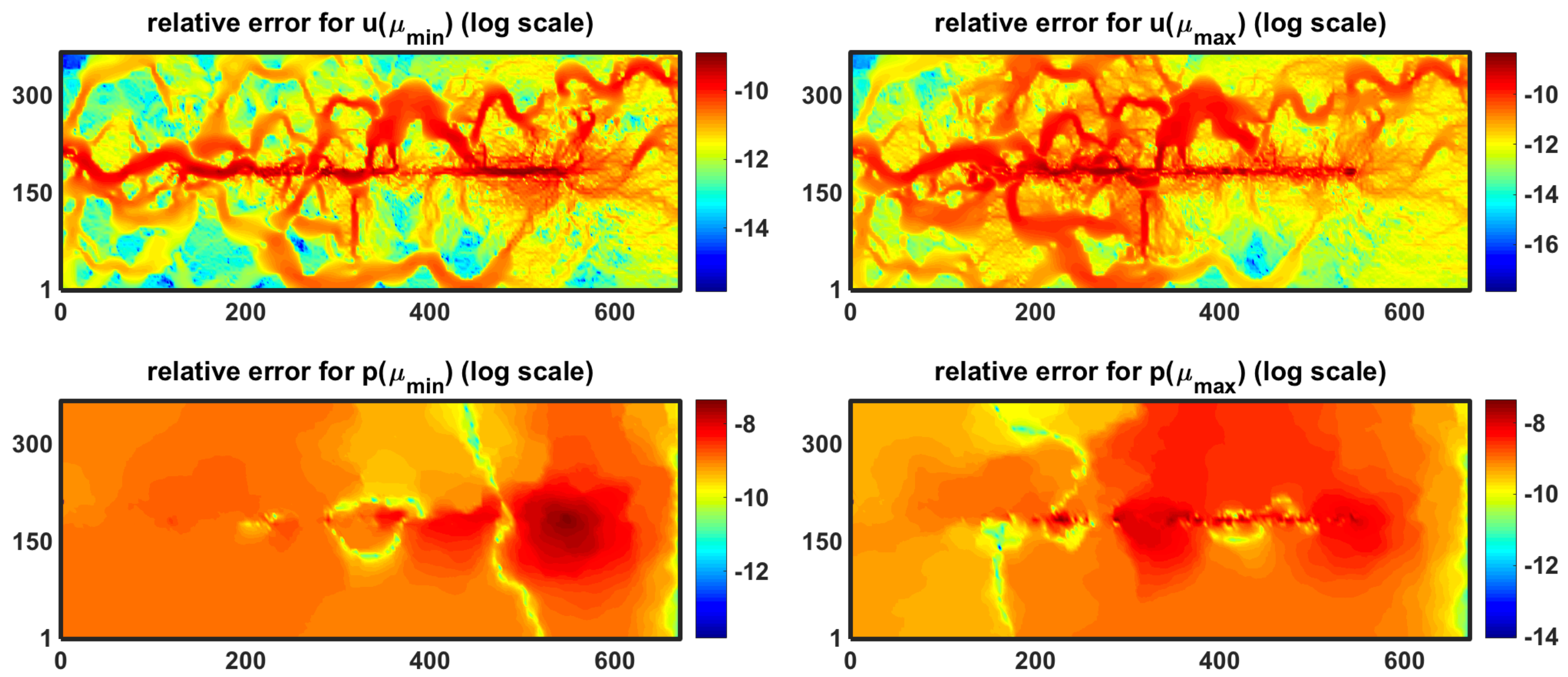}
 \caption{Fractured medium. Maps for the relative difference between PGD and corresponding FE solution for the norm of the velocity  and pressure (top and bottom rows respectively) evaluated at minimum and maximum value for the parameter (left and right columns respectively).}
 \label{fig:brinkman:error}
\end{figure*}

Here, the model of K\"onn\"o and Stenberg is extended including one additional parameter that controls the permeability of the fracture, occupying a subdomain $\Lambda$ in $\Omega$  . This parameter is denoted by $\mu$ and can be interpreted as a measure of the degree of saturation of some \emph{filling} inside the fracture. The influence of the input parameter $\mu$ ranging in $[0, \mu_{\text{max}}]$ in the resulting permeability distribution is given by the following expression:
\begin{equation}\label{eq:kxmu}
k(\bx,\mu) = \left\{
     \begin{array}{ll}
       10^{n(\bx) -\frac{\mu}{\mu_{\text{max}}} n(\bx)  + \mu} & \text{for } \bx \in \Lambda \\
       10^{n(\bx)} & \text{for } \bx \in \Omega \setminus \Lambda 
     \end{array}
   \right. 
\end{equation}
where $n(\bx) := \log_{10}\left(k_{\text{SPE10}}(\bx)\right)$. 
Note that for the extreme values of $\mu$,  $k(\bx,0)= k_{\text{SPE10}}(\bx)$ and $k(\bx,\mu_{\text{max}})=10^{\mu_{\text{max}}}$.

In order to implement the PGD, input data must be expressed in a separated form. The only term in \eqref{eq:brinkman} that it is not trivially separable is the second term of \eqref{eq:brinkman1}, involving the inverse of the permeability $k(\bx,\mu)$. Thus, the inverse of $k(\bx,\mu)$ as defined in \eqref{eq:kxmu} has to be expressed as a separated expression, that is in terms of functions that depend only on $n(\bx)$ and functions that depend only on $\mu$. Note that this is only needed for $ \bx \in \Lambda$.  The part of  $k(\bx,\mu)^{-1}$ which is not trivially separable is $10^{\frac{\mu}{\mu_{\text{max}}} n(\bx)}$. 
A SVD of a dense sampling of this function is used to separate it. The separated approximation is a description in terms of modes $\Theta^m(n(\bx))$ and $\phi^m(\mu)$, $m=1,\ldots,M$, and reads
\begin{equation*}
10^{\frac{\mu}{\mu_{\text{max}}} n(\bx)} \approx \sum_{m=1}^M \Theta^m(n(\bx))\, \phi^m(\mu)
\end{equation*}
Thus, for  $ \bx \in \Lambda$, 
\begin{equation*}
\eta k^{-1}(\bx,\mu) \approx \eta\,10^{-n(\bx)}   \left( \sum_{m=1}^M \Theta^m(n(\bx))\, \phi^m(\mu) \right)   10^{-\mu} .
\end{equation*}
Note that the separation is performed in terms of the variables $n(\bx)$ and $\mu$ (instead of $\bx$ and $\mu$). This simplifies the function to be separated and reduces the number of terms required to reach some prescribed accuracy. In this case, using $M=16$ terms provides a relative error smaller than $10^{-12}$ (for any value of $n(\bx)$ and $\mu$).  

The PGD solution is sought with the form defined in \eqref{eq:sep2} (same parameter functions for velocity and pressure). 
The evolution with the number of PGD terms of the error with respect to a standard FE solution is shown in Figure~\ref{fig:brinkman:convergence}.
It can be observed that the relative error is of order  $10^{-5}$ with only 40 PGD terms. 

Despite  being relatively simple problem (with only one scalar parameter $\mu$),  the convergence is faster than in other PGD solutions.
The resulting velocity and pressure fields for the extreme values of $\mu$ are shown in Figure~\ref{fig:brinkman:solution}, and the errors with respect to the corresponding finite element solutions are in Figure~\ref{fig:brinkman:error}. Note that for these two particular values, the errors are lower than for the worst case scenario depicted in the convergence curve of Figure~\ref{fig:brinkman:convergence}.


\section{Conclusions}
The analysis of the different forms for the parametric separation of Stokes problems reveals that the only viable option is having a unique parametric mode for each independent parameter, affecting all the velocity components and the pressure. This choice (denoted as case \#2) has less degrees of freedom than the alternative cases \#1 (one parametric mode for all velocity components and a different one for the pressure) and \#0 (different parametric modes for every velocity component and the pressure). 
The proposed formulation corresponding to case \#2 is the simplest alternative, guarantees incompressibility and is not affected by any stability concerns. 

The abundance of degrees of freedom may be a desirable feature of the PGD formulation because a richer functional space could compensate the computational  overhead for each mode with a reduced number of modes. However, in this case increasing the unknowns leads to unsolvable problems and therefore the alternatives \#0 and \#1 have to be discarded.

Moreover, the PGD compression algorithm based on a Least Squares projection can be performed for the three alternative parametric representations. The analysis of the three compressions demonstrates that the parametric structure of the solutions does not require the multiplicity of parametric modes: the number of modes required to reach some prescribed accuracy for case \#2 is almost the same as for the other two alternatives.
\section*{Acknowledgements}
This work has been supported by \emph{Ministerio de Econom\'{\i}a y Competitividad}, grant number DPI2014-51844-C2-2-R and  by \emph{Generalitat de Catalunya}, grant number 2014-SGR-1471.



\end{document}